\documentstyle[11pt,amssymb]{article}

\newtheorem{thm}{Theorem}[section]
\newtheorem{lem}[thm]{Lemma}
\newtheorem{prop}[thm]{Proposition}
\newtheorem{cor}[thm]{Corollary}
\newcommand{\Xast}{X^{\ast}}

\newcommand{\lip}{\mbox{\em lip}}
\newcommand{\norm}{\mbox{\em norm}}

\newcommand{\conv}{\mbox{\em conv}}
\newcommand{\dist}{\mbox{\em dist}}
\newcommand{\pf}{\noindent {\bf Proof.\ \ }}
\newcommand{\finpf}{\hfill{$\square$}\linespace}
\newcommand{\linespace}{\vspace{\baselineskip}  \noindent}
\newcommand{\con}{\mbox{conv}}

\newcommand{\nm}{\mbox{norm}}
\newcommand{\dn}{\mbox{dualnorm}}
\newcommand{\tn}{|\hskip-.13em | \hskip-.13em |}

\begin{document}

\title{\bf Locally Lipschitz Functions and Bornological Derivatives}
\author{J.M. Borwein$^{1}$ \\
Department of Combinatorics and Optimization \\
University of Waterloo, Waterloo, Ontario N2L 3G1 \\
\\
M. Fabian $^{2}$ \\
Department of Mathematics \\
Miami University, Oxford, Ohio 45056 \\
\\
J. Vanderwerff $^{3}$ \\
Department of Pure Mathematics \\
University of Waterloo, Waterloo, Ontario N2L 3G1}
\footnotetext[1]{Research supported in part by an NSERC (Canada)
operating grant.}
\footnotetext[2]{Research completed while visiting the Faculty of
Mathematics,
University of Waterloo.}
\footnotetext[3]{NSERC (Canada) postdoctoral fellow.}
\maketitle
\begin{abstract}
We study the relationships between Gateaux, Weak Hadamard and
Fr\'{e}chet
differentiability and their bornologies for Lipschitz and for convex
functions.
\end{abstract}
\vskip 2.0truecm
\noindent
{\it AMS Subject Classification.} Primary: 46A17, 46G05, 58C20.

\hskip 1.66truein Secondary: 46B20.

\medskip\noindent
{\it Key Words:} Gateaux differentiability, Fr\'echet
differentiability,
weak Hadamard differentiability, convex functions, Lipschitz functions,
bornologies, Dunford-Pettis property, Grothendieck property, Schur
property,
not containing $\ell_1$.
\newpage

\section{Introduction}
\setcounter{equation}{0}
We begin by introducing the definitions and notation that will be
used.  Unless
otherwise specified, $X$ is an infinite dimensional real Banach space
with norm
$\|\cdot\|$ and dual space $\Xast$.  A {\it bornology} ${\cal B}$ on
$X$ is a
family of bounded subsets of $X$ such that
$\cup \left\{ B:B \in {\cal B} \right\} =
X$.  We will focus on the following bornologies:  $G=$ \{singletons\},
$H=$
\{compact sets\}, $WH=$ \{weakly compact sets\} and $F=$ \{bounded
sets\}.
Observe that $G \subset H \subset WH \subset F$.

A function $f : X \rightarrow \Bbb{R}$ is called {\it ${\cal
B}$-differentiable}
at $x \in X$ if there is $\Lambda \in \Xast$
such that for each $B \in {\cal B}$,
\[ \frac{1}{t} \Bigl[ f(x+th) - f(x) -
\langle \Lambda, th\rangle \Bigr] \rightarrow 0~~ {\rm as}~~
t \downarrow 0 \]
uniformly for $h \in B$.  Let ${\cal F}$ denote a family of real-valued
locally
Lipschitz functions on $X$; we will usually consider locally Lipschitz
(loclip),
Lipschitz (lip), distance (dist), continuous convex (conv) and norms
(norm);
it is, of course, easy to check that continuous convex functions
are locally Lipschitz (\cite{Ph},
Proposition 1.6).
For two bornologies on a fixed Banach space $X$,
say $F$ and $G$ and a family of functions
${\cal F}$, we will write $F_{{\cal F}} = G_{{\cal F}}$,
if for every $f \in {\cal F}$ and every
$x \in X, f$ is $F$-differentiable at $x$ if
and only if $f$ is $G$-differentiable at $x$.  Since
$G_{\mbox{loclip}} = H_{\mbox{loclip}}$, we
will write $G$ and $H$ interchangeably.

In the paper \cite{BF}, it was shown that $H_{\con}
= F_{\con}$ if and only if
$X$ is infinite dimensional.  From this, one might be tempted to
believe that
various notions of differentiability
for convex functions coincide precisely when the
bornologies on the space coincide.  However, this is far from being the
case; for
example, according to (\cite{BF}, Theorem 2),
$WH_{\con} = F_{\con}$ if and only if
$X \not\supset \ell_1$.  In contrast to this, we will show in Section 2
that
differentiability notions coincide for Lipschitz functions precisely
when the
bornologies are the same (for the $H$, $WH$ and $F$ bornologies).
In the third section
we will study the relationship between
$WH$-differentiability and $H$-differentiability
for continuous convex functions. In particular, we show that if
$B_{X^*}$ is
$w^*$-sequentially compact, then $H_{\con} = WH_{\con}$ precisely when
$H = WH$.
However, there are spaces for which $H_{\con} = WH_{\con}$ and yet
$H \not= WH$. This leads to examples showing that one cannot always
extend
a convex function from a space to a superspace while preserving
$G$-differentiability at a prescribed point. Some characterizations of
the
Schur and Dunford-Pettis properties are also obtained in terms of
differentiabilty of continuous $w^*$-lower semicontinuous convex
functions
on the dual space.

\section{Lipschitz functions and bornologies}
\setcounter{equation}{0}
As mentioned in the introduction, there are spaces for which $WH \neq
F$ but
$WH_{\con} = F_{\con}$.   However, this is not
the case for Lipschitz functions.

\begin{thm}\label{thm2.1}
For a Banach space $X$, the following are equivalent.
\begin{description}
\item[(a)]
$X$ is reflexive.
\item[(b)]
$WH_{\lip} = F_{\lip}$.
\item[(c)]
$WH_{\dist} = F_{\dist}$.
\end{description}
\end{thm}

In order to prove this theorem, we will need a special type of sequence
in
nonreflexive Banach spaces.  Namely, we will say
$\left\{x_k \right\}^{\infty}_{k=1} \subset X$ is a {\it special
sequence} if there is an $\epsilon > 0$ such that
$\left\{ z_k \right\} \subset X$ has no weakly convergent
subsequence whenever $\|x_k - z_k \| < \epsilon$.

\medskip
\noindent
{\bf Remark}
(a) There are examples of sequences $\left\{ x_k \right\}$
such that $\left\{ x_k \right\}$
has no weakly convergent subsequence but $\left\{ x_k \right\}$ is not
special.

Indeed, let $X = \ell_1$ and consider
$y_{n,m} = e_n + \frac{1}{n}e_m$ for $m$, $n \in
\Bbb{N}, m > n$.  Let $\{x_k\}$ be any sequential arrangement of
$\left\{ y_{n,m} \right\}$.
It is not hard to verify $\left\{ x_k \right\}$ has the desired
properties.  Another example is $X = c_0$
and $y_{n,m} =  {\displaystyle \sum^{n}_{k=1}} e_k
+ {\displaystyle \sum^{n+m}_{k = n+1}} \frac{1}{n}e_k$.

(b) If $f$ is Lipschitz and $WH$-differentiable at $0$ (with
$f^{\prime}(0)=0)$ but $f$ is not Fr\'{e}chet
differentiable at 0, it is not hard to
construct a special sequence $\left\{ x_k \right\}$.

Indeed, because $f$ is not Fr\'{e}chet differentiable,
we can choose $\left\{ x_k
\right\}$ in the unit sphere $S_{X}$ of $X$ and $t_k \downarrow 0$
which satisfy
\[ \frac{|f(t_k x_k) - f(0)|}{t_k} \geq \epsilon \hspace{.5in} {\rm
for~~some}~~
\epsilon > 0. \]
Using the fact that $f$ is Lipschitz and
$WH$-differentiable at 0, one can easily
show that $\left\{ x_k \right\}$ is special.
\finpf

Part (b) of the above remark shows that in order to prove Theorem
\ref{thm2.1}, it is necessary
to show each nonreflexive Banach space has a special sequence.  On the
other hand,
part (a) shows that such sequences must be chosen carefully.

\begin{lem} \label{lem2.2}
Suppose $\left\{ x_n \right\} \subset X$ has
no weakly convergent subsequence.  Then some
subsequence of $\left\{x_n \right\}$ is a special sequence.
\end{lem}

\pf
If some subsequence of $\left\{ x_n \right\}$ is special, then
there is nothing more to
do.  So we will suppose this is not so and arrive at a contradiction by
producing a
weakly convergent subsequence of $\left\{ x_n \right\}$.

Given $\epsilon = 1$, by our supposition, we choose $N_1 \subset
\Bbb{N}$ and $\left\{
z_{1,i} \right\}_{i \in N_1}$ such that
\[ \| x_i - z_{1,i} \| < 1 ~~ \enskip ~~ {\rm for} ~~ \enskip ~~ i \in
N_1 ~~\enskip~~
{\rm and}~~\enskip~~  w\!{\rm-}\!\!\lim_{i \in N_1}  z_{1,i} = z_1 . \]
Supposing $N_{k-1}$ has been chosen, we choose $N_k \subset N_{k-1}$
and $\left\{ z_{k,i}
\right\}_{i \in N_k} \subset X$ satisfying
\begin{equation}\label{eqn2.1}
\| x_i - z_{k,i} \| < \frac{1}{k} ~~\enskip~~ {\rm for}~~ i \in N_k
~~\enskip~~ {\rm and}
~~\enskip~~ w\!{\rm-}\!\!\lim_{i \in N_k} z_{k,i} = z_k .
\end{equation}
In this manner we construct $\left\{ z_{k,i} \right\}_{i \in N_k}$ and
$N_k$
for all $k \in \Bbb{N}$.

Notice that $z_n - z_m = {\displaystyle w\!{\rm-}\!\!\lim_{i \in N_n}}
(z_{n,i} - z_{m,i})$ for $n > m$.
Thus by the $w$-lower semicontinuity of $\|\cdot\|$ and (\ref{eqn2.1})
we obtain
\[ \| z_n - z_m \| \leq \liminf_{i \in N_n} \|z_{n,i} - z_{m,i} \| \leq
\liminf_{i \in N_n} (\|
z_{n,i} - x_i \| + \| x_i - z_{m,i} \|) \leq \frac{1}{n} + \frac{1}{m}
\leq \frac{2}{m}. \]
Thus $z_n$ converges in norm to some $z_{\infty} \in X$.

Now for each $n \in \Bbb{N}$ choose integers $i_n \in N_n$ with $i_n >
n$.
We will show $x_{i_{n}} \stackrel{w}{\rightarrow}
z_{\infty}$.  So let $\Lambda \in B_{X^{\ast}}$ and $\epsilon > 0$ be
given.  We select an $n_0 \in
\Bbb{N}$ which satisfies
\begin{equation}\label{eqn2.2}
\frac{1}{n_0} < \frac{\epsilon}{3} ~~\enskip~~ {\rm and} ~~\enskip~~ \|
z_m -
z_{\infty} \| < \frac{\epsilon}{3} ~~\enskip~~ {\rm for} ~~\enskip~~ m
\geq n_0.
\end{equation}
Because $z_{n_{0},i} \stackrel{w}{\rightarrow} z_{n_{0}}$,
we can select $m_0$ so that
\begin{equation}\label{eqn2.3}
\bigl| \langle \Lambda, z_{n_{0},i} - z_{n_{0}} \rangle
\bigr|  < \frac{\epsilon}{3} ~~ \enskip~~
{\rm for~all} ~~\enskip~~ i \geq m_0.
\end{equation}
For $m \geq \max \left\{ n_0, m_0 \right\}$, we have
\[ \begin{array}{ll}
\bigl| \langle \Lambda, x_{i_{m}} - z_{\infty} \rangle \bigr|
& \leq \bigl| \langle \Lambda,
x_{i_{m}} - z_{n_{0},i_{m}} \rangle \bigr| +
\bigl| \langle \Lambda, z_{n_{0},i_{m}} - z_{n_{0}}
\rangle \bigr| + \bigl| \langle \Lambda,
z_{n_{0}} - z_{\infty} \rangle \bigr| \\
& < \| x_{i_{m}} - z_{n_{0},i_{m}} \| + \frac{\epsilon}{3} + \|
z_{n_{0}} - z_{\infty} \| \hspace{.25in}
\left[ {\rm by~(\ref{eqn2.3})~since}~i_m > m \geq m_0 \right] \\
& < \frac{1}{n_0} + \frac{\epsilon}{3} + \frac{\epsilon}{3} <
\epsilon.
\hspace{2.1in} {\rm [by~
(\ref{eqn2.2})~ and ~(\ref{eqn2.1})]}
\end{array} \]
Therefore $x_{i_{n}} \stackrel{w}{\rightarrow} z_{\infty}.$
\finpf

\noindent
{\bf Proof of Theorem 2.1.}  Notice that (a) $\Longrightarrow$ (b)
$\Longrightarrow$ (c) is trivial.  It
remains to prove (c) $\Longrightarrow$ (a).
Suppose $X$ is not reflexive, hence $B_X$ is not weakly
compact and so there exists $\left\{ x_n \right\} \subset S_X$
with no weakly convergent subsequence.  By
Lemma 2.2 there is a subsequence, again denoted by $\left\{ x_n
\right\}$,
and a $\Delta \in (0,1)$ such that
$ \left\{ z_n \right\} \subset X$ has no weakly convergent
subsequence whenever $\| z_n - x_n \| < \Delta$.  By
passing to another subsequence if necessary
we may assume $\| x_n - x_m \| > \delta$
for all $n \neq m$, with some $0<\delta<1$.

For $n = 1,2, \ldots,$ let $B_n = \{x \in X : \|x - 4^{-n}x_n\| \le
\delta \Delta 4^{-n-1} \}$ and put $C = X \backslash
{\displaystyle \cup^{\infty}_{n=1}} B_n$.  Because $4^{-m} + \delta
\Delta 4^{-m-1} < 4^{-n} - \delta
\Delta 4^{-n-1}$ for $m > n$, we have that $B_n \cap B_m = \emptyset$
whenever $n \neq m$.  For $x \in X$,
let $f(x)$ be the distance of $x$ from $C$.  Thus $f$ is a Lipschitz
function on $X$ with $f(0) = 0$.  We
will check that $f$ is $WH$-differentiable at 0 but not
$F$-differentiable at 0.

Let us first observe that $f$ is $G$-differentiable at 0.  So fix any
$h \in X$ with $\|h\|=1$.  Then
$[0,+\infty)h$ meets at most one ball $B_n$.  In fact assume $t_m, t_n
> 0$ are such that $\|t_{i}h
- 4^{-i}x_i \| < \delta \Delta 4^{-i-1}$ for $i = n,m$.  Then
$|4^{i}t_i - 1 | < \frac{\delta \Delta}{4}$
for $i = n,m$ and
\[ \begin{array}{ll}
\| x_n - x_m \| & \leq \| x_n - 4^n t_n h\|
+ \| 4^n t_n h - 4^m t_m h \| + \| 4^m t_m h - x_m \| \\
& < \frac{\delta \Delta}{4} + \frac{2 \delta \Delta}{4} + \frac{\delta
\Delta}{4} = \delta \Delta < \delta.
\end{array} \]
This means that $n=m$.  It thus follows that for $t>0$ small enough,
we have $f(th)=0$.  Therefore $f$ is
$G$-differentiable at 0, with $f^{\prime}(0) = 0$.  Let us further
check that $f$ is not
$F$-differentiable at 0.  Indeed,
\[ \frac{f(4^{-n}x_n)}{\| 4^{-n}x_n\|} = \frac{\delta \Delta}{4}
~~\enskip~~ {\rm for~all}~~n, \]
while $\|4^{-n}x_n \| \rightarrow 0$.

Finally assume that $f$ is not $WH$-differentiable at 0.  Then there
are a weakly compact set $K \subset
B_X, \epsilon > 0$, and sequences $\left\{ k_m \right\} \subset K,
t_m \downarrow 0$ such that
\[ \frac{f(t_m k_m)}{t_m} > \epsilon ~~\enskip~~ {\rm for~all}~~ m \in
\Bbb{N}. \]
Hence, as $f$ is $1$-Lipschitz, we have $\inf \| k_n \| \ge \epsilon>
0$.
Further,
because $f(t_m k_m) > 0$, there are $n_m \in \Bbb{N}$ such that
\[ \| t_m k_m - 4^{-n_{m}} x_{n_{m}} \|
< \Delta \delta 4^{-m_{n}-1}, ~~ m = 1,2, \ldots\ . \]
Consequently,
\begin{equation}\label{eqn2.4}
\|4^{n_{m}}t_m k_m - x_{n_{m}} \|< \frac{\Delta \delta}{4} < \Delta
~~\mbox{and}~~ |4^{n_{m}}
t_m \| k_m \| - 1 | < \frac{\Delta \delta}{4} .
\end{equation}
Because $\left\{ x_n \right\}$ is a special sequence with $\Delta$, the
first inequality in (\ref{eqn2.4})
says that $\left\{ 4^{n_{m}} t_m k_m \right\}$
does not have a weakly convergent subsequence.  However the
second inequality in (\ref{eqn2.4}) together with $\inf \| k_n \| > 0$
ensures
that  $4^{n_m}t_m$ is bounded and so, since $\{k_m\}$ is weakly
compact,
$4^{n_m}t_m k_m$ has a weakly convergent subsequence,
a contradiction.  This proves $f$ is $WH$-differentiable at 0.
{\hfill{$\square$}\linespace}

Recall that a Banach space has the {\it Schur property} if $H = WH$,
that is, weakly convergent sequences
are norm convergent.

\begin{thm}\label{thm2.3}
For a Banach space $X$, the following are equivalent.
\begin{description}
\item[(a)]
$X$ has the Schur property.
\item[(b)]
$H_{\lip} = WH_{\lip}.$
\item[(c)]
$H_{\dist} = WH_{\dist}$.
\end{description}
\end{thm}

\pf
It is clear that (a) $\Longrightarrow$ (b) $\Longrightarrow$ (c), thus
we prove (c) $\Longrightarrow$
(a).  Suppose $X$ is not Schur and choose $\left\{ x_n \right\} \subset
S_X$ such that $x_n
\stackrel{w}{\rightarrow} 0$ but $\| x_n \| \not\rightarrow 0$.  Since
$\left\{ x_n \right\}$ is not
relatively norm compact, we may assume by passing to a subsequence if
necessary that $\| x_i - x_j \| >
\delta$ for some $\delta \in (0,1)$ whenever $i \neq j$.

As in the proof of Theorem 2.1, let $B_n = \{x \in X: \|x - 4^{-n}x_n\|
\le \delta 4^{-n-1}\}$, $C = X \backslash {\displaystyle
\cup^{\infty}_{n = 1}} B_n$ and let $f(x) = d(x,C)$. Now $f(0)=0$ and
the argument of Theorem \ref{thm2.1}, shows that $f$ is
$G$-differentiable at 0 with $f^{\prime}(0)=0$.  However,
\[ \frac{f(4^{-n}x_n)}{4^{-n}} = \frac{\delta}{4} ~~\enskip~~ {\rm
for~all} ~~ n \in \Bbb{N}. \]
Since $\left\{ x_n \right\} \cup \left\{ 0 \right\}$ is weakly compact,
it follows that $f$ is not
$WH$-differentiable at 0.
\finpf

\noindent
{\bf Remark.} Using the technique from the proof of Theorem 2.1,
one can also prove
the following statement. If a nonreflexive
Banach space $X$ admits a Lipschitzian $C^k$-smooth bump function,
then it admits a Lipschitz function which is $C^k$-smooth on $X
\backslash \left\{ 0 \right\},
WH$-differentiable at 0, but not $F$-differentiable at 0.
A corresponding remark holds for non-Schur spaces.

\section{Differentiability properties of convex functions}
\setcounter{equation}{0}
We begin by summarizing some known results.  First recall that a Banach
space
$X$ has the {\it Dunford-Pettis property} if
$\langle x^{\ast}_{n}, x_{n}\rangle \rightarrow 0$
whenever $x_{n} \stackrel{w}{\rightarrow} 0$ and $x^{\ast}_{n}
\stackrel{w}{\rightarrow} 0$.  For notational
purposes we will say $X$ has the
$DP^{\ast}$ if $\langle x^{\ast}_{n},x_n \rangle
\rightarrow 0$ whenever $x^{\ast}_{n}
\stackrel{w^{\ast}}{\rightarrow} 0$ and $x_n \stackrel{w}{\rightarrow}
0$; see
(\cite{DU}, p. 177) for more on the Dunford-Pettis property.  Note that
a {\it
completely continuous} operator takes weakly convergent sequences to
norm convergent
sequences. The proof of the next result is essentially in {\cite{BF}}.

\begin{thm}(\cite{BF}) \label{thm3.1}
\begin{enumerate}
\item[(a)]
$X$ does not contain a copy of $\ell_1$ if and only if $WH_{\conv} =
F_{\conv}$
if and only if $WH_{\norm} = F_{\norm}$
if and only if each completely continuous linear $T : X \rightarrow
c_0$ is
compact.
\item[(b)]
$X$ has the $DP^{\ast}$ if and only if $H_{\conv} = WH_{\conv}$ if and
only if
$H_{\norm} = WH_{\norm}$ if and only if
each continuous linear $T : X \rightarrow c_0$ is completely
continuous.
\item[(c)]
$X$ is finite dimensional if and only if $G_{\conv} = F_{\conv}$ if and
only if
$G_{\norm} = F_{\norm}$ if and only if
each continuous linear $T : X \rightarrow c_0$ is compact.
\end{enumerate}
\end{thm}

\pf Let us mention that (a) is contained in (\cite{BF}, Theorem 2)
and (c) is from
(\cite{BF}, Theorem 1). Whereas (b) can be obtained by following the
proofs of (\cite{BF}, Proposition 1 and Theorem 1).
\finpf

If $WH_{\con} \neq G_{\con}$, for example, we can be
somewhat more precise.

\begin{prop}\label{prop3.2}
Suppose $WH_{\conv} \neq G_{\conv}$ on $X$.
Then there is a norm $\tn \cdot \tn$ on
$X$ such that $\tn \cdot \tn$ is not $WH$-differentiable
at $x_0 \neq 0$ but $\tn \cdot
\tn^{\ast}$ is strictly convex at
$\Lambda_0 \in X^{\ast}\backslash\{0\}$
satisfying $\langle \Lambda_0,x_{0}\rangle
=\tn x_0\tn\, \tn \Lambda_0\tn $.

\end{prop}

\pf
Following the techniques of \cite{BF}, one obtains a norm $\| \cdot \|$
on $X$
such that $\| \cdot \|$ is $G$-differentiable at $x_0 \neq 0$ but $\|
\cdot \|$
is not $WH$-differentiable at $x_{0}$.  Now define $\tn \cdot \tn$ on
$X$ by
\[ \tn x \tn = (\|x\|^2 + d^2 (x,\Bbb{R} x_0)
)^{\frac{1}{2}} \]
Clearly $d^2 ( \cdot , \Bbb{R} x_0 )$ is $F$-differentiable at
$x_0$ and so it follows that $\tn \cdot \tn$ is $G$-differentiable at
$x_0$
but $\tn \cdot \tn$ is not $WH$-differentiable at $x_0$ because $\|
\cdot \|$
is not.  Suppose now that $\left\{ x_n \right\}$ satisfies
\begin{equation}\label{eqn3.1}
2 \tn x_n \tn^2 + 2\tn x_0 \tn^2 - \tn x_n + x_0 \tn^2 \rightarrow 0.
\end{equation}
Then by convexity one obtains
\[ \| x_n \| \rightarrow \| x_0 \| ,\  d^2 (x_n , \Bbb{R} x_0 )
\rightarrow d^2
(x_0 , \Bbb{R} x_0 ) = 0.
\]
>From this one easily sees that $\| x_n - x_0 \| \rightarrow 0.$

Now take $\Lambda_0 \in X^*$ such that $\tn \Lambda_0\tn =1$ and
$\langle \Lambda_0, x_0 \rangle = \tn x_0 \tn$.
We show that $\tn \cdot \tn$ is
strictly convex at $\Lambda_0$.
Suppose that $\tn x^*\tn  = 1$ and $\tn x^* + \Lambda_0\tn = 2$, then
choose $\left\{ x_n \right\}$ with $\tn x_n \tn = \tn x_0 \tn$ so that
\begin{equation}\label{eqn3.2}
\langle x^* + \Lambda_0 , x_n \rangle \rightarrow 2 \tn x_0 \tn.
\end{equation}
Consequently $\langle \Lambda_0 , x_n \rangle
\rightarrow \tn x_0 \tn$ and thus
$\langle \Lambda_0 , x_n + x_0 \rangle \rightarrow 2\tn x_0 \tn$;
$\tn x_n + x_0 \tn \to 2\tn x_0 \tn$.  But then
$\left\{ x_n \right\}$ satisfies (\ref{eqn3.1}) and so $\| x_n - x_0 \|
\rightarrow 0$.  This with (\ref{eqn3.2}) shows
that $\langle x^* , x_0 \rangle =
\tn x_0 \tn$.  Because $\tn \cdot \tn$ is $G$-differentiable at $x_0$,
we
conclude that $x^* = \Lambda_0$.  This proves the strict convexity of
$\tn \cdot \tn$ at $\Lambda_0$.
\finpf

One can also formulate similar  statements (and proofs) for the
cases $G_{\con} \neq F_{\con}$ and
$WH_{\con} \neq F_{\con}$.

We now turn our attention to spaces for which $WH_{\con} = H_{\con}$.
Let us
recall that a Banach space $X$ has the {\it Grothendieck property} if
$w^{\ast}$-convergent sequences in $X^{\ast}$ are weakly convergent;
see
(\cite{DU}, p. 179).  The following corollary is an immediate
consequence of
Theorem \ref{thm3.1} (b).

\begin{cor} \label{cor3.3}
If $X$ has the Dunford-Pettis property and the Grothendieck property,
then
$WH_{\conv} = H_{\conv}$.
\end{cor}

In particular, note that $\ell_{\infty}$ has the Grothendieck property
(cf
\cite{D}, p.103) and the Dunford-Pettis property (cf \cite{DU}, p.
177).  Thus,
unlike the case for Lipschitz functions, one can have $WH_{\con} =
H_{\con}$
for non-Schur spaces.  It will follow from the next result
that these non-Schur spaces
must be quite large, though.

\begin{thm} \label{thm3.4}
For a Banach space $X$, the following are equivalent.
\begin{itemize}
\item[(i)] $X$ has the $DP^{\ast}$
\item[(ii)] $H_{\conv} = WH_{\conv}$
\item[(iii)] If $B_{Y^{\ast}}$ is $w^{\ast}$-sequentially compact,
then any continuous linear $T : X \rightarrow Y$ is completely
continuous.
\end{itemize}
\end{thm}

\pf
By Theorem \ref{thm3.1}(b) we know that (i) and (ii) are equivalent and
that
(iii) implies (i).  We will show (i) implies (iii) by contraposition.
Suppose
(iii) fails, that is, there is an operator $T : X \rightarrow Y$ which
is not
completely continuous for some $Y$ with $B_{Y^{\ast}}$
$w^{\ast}$-sequentially
compact.  Hence we choose $\left\{ x_n \right\} \subset X$ such that
$x_n
\stackrel{w}{\rightarrow} 0$ but $\| T x_n \| \not\rightarrow 0$.
Because $T
x_n \stackrel{w}{\rightarrow} 0$, we know that $\left\{ T x_n \right\}$
is not
relatively norm compact.  Hence letting $E_n = {\rm span} \left\{ y_k :
k \leq n
\right\}$ with $y_k = T x_k$ we know there is an $\epsilon > 0$ such
that
${\displaystyle \sup^{}_{k}} \,d (y_k , E_n ) > \epsilon$ for each
$n$.

By passing to a subsequence, if necessary, we assume $d(y_n , E_{n-1})
>
\epsilon$ for each $n$.  Now choose $\Lambda_n \in B_{Y^{\ast}}$ such
that
$\langle \Lambda_n , x \rangle = 0$ for all $x \in E_{n-1}$
and $\langle \Lambda_n , x_n \rangle
\geq \epsilon$.  Because $B_{Y^{\ast}}$ is $w^{\ast}$-sequentially
compact,
there is a subsequence $\Lambda_{n_{k}}$ such that $\Lambda_{n_{k}}
\stackrel{w^{\ast}}{\rightarrow} \Lambda \in B_{Y^{\ast}}$.
Observe that $\langle \Lambda_{n}, y_{k}\rangle =
0$ for $n > k$ and consequently $\langle \Lambda, y_k\rangle = 0$
for all $k$.  Now let
$z^{\ast}_{k} = T^{\ast} (\Lambda_{n_{k}} - \Lambda)$
and $z_k = x_{n_{k}}$.  Certainly
$z^{\ast}_{k} \stackrel{w^{\ast}}{\rightarrow} 0$
and $z_k \stackrel{w}{\rightarrow}
0$ while $\langle z^{\ast}_{k},z_k \rangle = \langle \Lambda_{n_{k}} -
\Lambda, Tx_{n_{k}}\rangle = \langle \Lambda_{n_{k}} -
\Lambda, y_{n_{k}}\rangle \geq \epsilon$
for all $k$.  This shows that $X$ fails the
$DP^{\ast}$.
\finpf

\begin{cor} \label{cor3.5}
If $X$ has a $w^{\ast}$-sequentially compact dual ball or,
more generally, if every separable
subspace of $X$ is a subspace of a complemented subspace with
$w^{\ast}$-sequentially compact dual ball, then the following are
equivalent.
\begin{itemize}
\item[(a)] $WH_{\conv} = H_{\conv}$.
\item[(b)] $X$ has the Schur property.
\end{itemize}
\end{cor}

\pf
Note that (b) $\Longrightarrow$ (a) is always true, so we show (a)
$\Longrightarrow$ (b).  If $B_{X^{\ast}}$ is $w^{\ast}$-sequentially
compact
and $X$ is not Schur then $I : X \rightarrow X$ is not completely
continuous
and Theorem 3.4 applies.
More generally, suppose $x_n \stackrel{w}{\rightarrow} 0$ but $\| x_n
\|
\not\rightarrow 0$ and ${\overline {\rm span}}
\left\{ x_n \right\} \subset Y$ with
$B_{Y^{\ast}}$ $w^{\ast}$-sequentially compact.  If there is a
projection $P :
X \rightarrow Y$, then $P$ is not completely continuous since $P|_{Y}$
is the identity on $Y$.
\finpf

We can say more in the case that $X$ is weakly countably determined
(WCD);
see \cite{M} and Chapter VI of \cite{DGZ} for the definition and
further
properties of WCD spaces.
\begin{cor} \label{cor3.6}
For a Banach space $X$, the following are equivalent.
\begin{itemize}
\item[(a)] $X$ is WCD and $WH_{\conv} =
H_{\conv}$
\item[(b)] $X$ is separable and has the Schur property.
\end{itemize}
\end{cor}

\pf
It is obvious that (b) $\Longrightarrow$ (a), so we show (a)
$\Longrightarrow$
(b).  First, since $B_{X^{\ast}}$ is $w^{\ast}$-sequentially compact
(see e.g. \cite{M}, Corollary 4.9 and \cite{LP}, Theorem 11),
it follows from Corollary \ref{cor3.5} that $X$ has the Schur
property.  But $WCD$ Schur spaces are separable (see e.g. \cite{M},
Theorem
4.3).
\finpf

\noindent
{\bf Remark}
\begin{enumerate}
\item[(a)] Corollary \ref{cor3.5} is satisfied, for instance, by $GDS$
spaces (see \cite{LP}, Theorem 11) and spaces
with countably norming $M$-basis (see \cite{Pl}, Lemma 1).
Notice that $\ell_1 (\Gamma)$ has a
countably norming $M$-basis for any $\Gamma$,
thus spaces with countably norming $M$-bases and the Schur property
need not be
separable.
\item[(b)]  If $X^{\ast}$ satisfies $WH_{\con} = F_{\con}$, then $X$
also does
(because $L_1 \subset X^{\ast}$ if $\ell_1 \subset X$ (see \cite{Dul}
Proposition 4.2)) but not conversely ($c_0$ and
$\ell_1$); cf. Theorem \ref{thm3.1}(a).
\item[(c)]  Let $X$ be a space such that $X$
is Schur but $X^{\ast}$ does not have
the Dunford-Pettis property (cf. \cite{DU}, p 178).
Then $X$ satisfies $H_{\con} =
WH_{\con}$ but $X^{\ast}$ does not satisfy $H_{\con} = WH_{\con}$.
\item[(d)]  There are spaces with the $DP^{\ast}$ that
are neither Schur nor have the
Grothendieck property; for example $\ell_1 \times \ell_{\infty}$.
\item[(e)]  It is well-known that $\ell_{\infty}$ has $\ell_2$ as a
quotient
(\cite{LT}, p. 111).
Thus quotients of spaces with the
$DP^{\ast}$ need not have the $DP^{\ast}$. It is clear that superspaces
of
spaces with the $DP^*$ need not have the $DP^*$; the example $c_0
\subset \ell_\infty$ shows that subspaces need not inherit the $DP^*$.

\item[(f)] Haydon (\cite{H}) has constructed a nonreflexive
Grothendieck $C(K)$
space which does not contain $\ell_\infty$. Using the continuum
hypothesis,
Talagrand (\cite{T}) constructed a nonreflexive Grothendieck $C(K)$
space $X$
such that $\ell_\infty$ is neither a subspace nor a quotient of $X$.
Since
$C(K)$ spaces have the Dunford-Pettis property (see \cite{D}, p. 113),
both
these spaces have the $DP^*$.
\end{enumerate}

As a byproduct of Corollaries \ref{cor3.3} and \ref{cor3.5}
we obtain the following
example which is related to results from (\cite{Z}).

\smallskip
\noindent
{\bf Example.}  Let $X$ be a space with the Grothendieck and
Dunford-Pettis
properties such that $X$ is not Schur (e.g. $\ell_{\infty}$).  Then
there is a
separable subspace $Y$ (e.g. $c_0$) of $X$ and a continuous convex
function $f$ on $Y$ such that
$f$ is $G$-differentiable at 0 (as a function on $Y$), but no
continuous convex
extension of $f$ to $X$ is $G$-differentiable at 0 (as a function on
$X$);
there also exist $y_0
\in Y\backslash\{0\}$ and an equivalent norm $\|\cdot\|$ on $Y$ whose
dual
norm is strictly convex but no extension of $\|\cdot\|$ to $X$ is
G-differentiable at $y_0$.

\smallskip
\pf
Let $Y$ be a separable non-Schur subspace of $X$.  By Corollary
\ref{cor3.5},
there is a continuous convex function $f$ on $Y$
which is $G$-differentiable at 0, but is not $WH$-differentiable at 0.
Since any extension $\tilde{f}$ of $f$ also fails to be
$WH$-differentiable at 0, it follows that $\tilde{f}$ is not
$G$-differentiable at 0
because $X$ has the $DP^{\ast}$. Because $Y$ fails the $DP^*$, there is
a
sequence $\{ \Lambda_n \} \subset X^*$ such that $\Lambda_n$ converges
$w^*$ but not Mackey to $0$. By the proof of (\cite{BF}, Theorem 3),
there
is a norm $\|\cdot\|$ on $Y$ whose dual is strictly convex that fails
to be $WH$-differentiable at some $y_0 \in Y\backslash \{0\}$; as
above, no
extension of $\|\cdot\|$ to $X$ can be $G$-differentiable at $y_0$.
\finpf

We close this note by relating the Schur and Dunford-Pettis properties
to some
notions of differentiability for dual functions.

\begin{thm} \label{thm3.7}
For a Banach space $X$, the following are equivalent.
\begin{itemize}
\item[(a)] $X$ has the Schur property.
\item[(b)] $G$-differentiability and $F$-differentiability coincide for
$w^{\ast}$-$\ell sc$ continuous convex functions on $X^{\ast}$.
\item[(c)] $G$-differentiability and $F$-differentiability coincide for
dual norms on
$X^{\ast}$.
\item[(d)] $H_{\lip} = WH_{\lip}$.
\end{itemize}
\end{thm}

\pf
Of course (a) and (d) are equivalent according to Theorem
\ref{thm2.3}.

(a) $\Longrightarrow$ (b):  Suppose (b) does not hold.
Then for some continuous convex
$w^{\ast}$-$\ell sc$ $f$ on $X^{\ast}$,
there exists $\Lambda_{0} \in
X^{\ast}$ such that $f$ is $G$-differentiable at $\Lambda_{0}$ but $f$
is not
$F$-differentiable at $\Lambda_{0}$.
Let $f^{\prime} (\Lambda_{0}) = x^{\ast\ast} \in
X^{\ast\ast}$.  We also choose $\delta > 0$ and $K > 0$ such that for
$x_1^*, x_2^* \in B(\Lambda_{0},\delta)$ we have $|f(x_1^*) - f(x_2^*)|
\leq K\|x_1^* - x_2^*\|$ (since $f$ is locally Lipschitz).
Because $f$ is not $F$-differentiable at
$\Lambda_{0}$, there exist $t_n \downarrow 0,
t_n < \frac{\delta}{2}, \Lambda_{n} \in
S_{X^{\ast}}$ and $\epsilon > 0$ such that
\begin{equation} \label{eqn3.3}
f(\Lambda_{0} + t_n \Lambda_{n}) - f(\Lambda_{0}) -
\langle x^{\ast\ast},t_n \Lambda_{n} \rangle \geq \epsilon t_n.
\end{equation}
Because $f$ is convex and $w^{\ast}$-$\ell sc$, using the separation
theorem
we can choose $x_n \in X$ satisfying
\begin{equation}\label{eqn3.4}
\langle x_n,x^* \rangle \leq f(\Lambda_{0} + t_n \Lambda_{n} + x^*)
- f(\Lambda_{0} + t_n
\Lambda_{n}) + \frac{\epsilon t_n}{2} ~~ \mbox{for~all}~~ x^* \in
X^{\ast};
\end{equation}
Putting $x^* = -t_n \Lambda_n$ in
(\ref{eqn3.4}) and using (\ref{eqn3.3}) one obtains
\[ \begin{array}{lcl}
\langle x_n , t_n \Lambda_{n} \rangle & \geq
& f(\Lambda_{0} + t_n \Lambda_{n}) - f(\Lambda_{0}) -
\frac{\epsilon t_n}{2} \\
& \geq  & \langle x^{\ast\ast} ,
t_n \Lambda_{n} \rangle + \frac{\epsilon t_n}{2}.
\end{array} \]
And hence, $\|x_n - x^{\ast\ast} \| \geq \frac{\epsilon}{2}$ for all
$n$.

Let $\eta > 0$ and fix $x^* \in S_{X^*}$. Since $f$ is
$G$-differentiable
at $\Lambda_0$, there is a $0 < t_0 < \frac{\delta}{2}$ such that for
$|t| \le t_0$ we have
\begin{equation}\label{eqn3.5}
\langle x^{**}, t x^* \rangle - f(\Lambda_0 + t x^*) + f(\Lambda_0) \ge
- \frac{\eta}{2} t_0.
\end{equation}
Using (\ref{eqn3.4}) with the fact that $f$ has Lipschitz constant $K$
on
$B(\Lambda_0, \delta)$, for $|t| \le t_0$ we obtain
\[ \begin{array}{lcl}
\langle x_n, t x^* \rangle & \le & f(\Lambda_0 + t_n \Lambda_n + tx^*)
-
f(\Lambda_0 + t_n \Lambda_n) + \frac{\epsilon t_n}{2} \\
& \le & f(\Lambda_0 + t x^*) - f(\Lambda_0) + \frac{\epsilon t_n}{2} +
2K t_n.
\end{array} \]
Choosing $n_0$ so large that $\frac{\epsilon t_n}{2} + 2K t_n <
\frac{\eta}{2}
t_0$ for $n \ge n_0$, the above inequality yields
\begin{equation}\label{eqn3.6}
f(\Lambda_0 + tx^*) - f(\Lambda_0) - \langle x_n, t x^* \rangle \ge
-\frac{\eta}{2}t_0 ~~\mbox{for}~~ n\ge n_0,~|t|\le t_0.
\end{equation}
Adding (\ref{eqn3.5}) and (\ref{eqn3.6}) results in
\[
\langle x^{**} - x_n, t x^* \rangle \ge -\eta t_0 ~~\mbox{for}~~ n\ge
n_0, ~
|t| \le t_0.
\]
Hence $|\langle x^{**} - x_n, x^* \rangle| \le \eta$ for $n  \ge n_0$.
This shows that
$x_n \stackrel{w^{\ast}}{\rightarrow} x^{\ast\ast}$.
Combining this with the fact that
$\|x_n - x^{\ast\ast} \| \not\rightarrow 0$ shown above,
we conclude that for some $\delta >
0$ and some subsequence we have $\| x_{n_{i}} - x_{n_{i+1}} \| >
\delta$ for all
$i$.  However $x_{n_{i}} - x_{n_{i+1}}\stackrel{w}{\rightarrow} 0$
(in $X$) because
$x_{n_{i}} - x_{n_{i+1}} \stackrel{w^{\ast}}{\rightarrow} 0$ (in
$X^{\ast\ast}$).  This shows that $X$ is not Schur.

Since (b) $\Longrightarrow$ (c) is obvious, we show
that (c) $\Longrightarrow$ (a).
Write $X = Y \times \Bbb{R}$ and suppose that
$X$ is not Schur.  Then we can choose
$\left\{y_n \right\} \subset Y$ such that $y_n
\stackrel{w}{\rightarrow} 0$ but
$\|y_n \| =1$ for all $n$.  Let $\{\gamma_n\} \subset (\frac{1}{2}, 1)$
be such that $\gamma_n \uparrow 1$ and define $\tn \cdot \tn$
on $X^{\ast} = Y^{\ast} \times \Bbb{R}$ by
\[ \tn (\Lambda , t) \tn = \sup \bigl\{ |\langle \Lambda ,
 y_n \rangle + \gamma_n t| \bigr\}
\lor \frac{1}{2}(\|\Lambda\| + | t|).
\]
This norm is dual since it is a supremum of $w^{\ast}$-$\ell sc$
functions and the proof of (\cite{BF}, Theorem 1) shows that
$\tn \cdot \tn$ is Gateaux but not Fr\'echet differentiable at
$(0,1)$.
\finpf

If $X$ is not Schur, then the previous theorem ensures the existence of
a
$w^*$-$\ell sc$ convex continuous function on $X^*$ which is
G-differentiable
but not F-differentiable at some point. The following remark shows that
we can
be more precise if $X \not\supset \ell_1$.

\medskip
\noindent
{\bf Remark.}  If $X \not\supset \ell_1$ and $X$ is not reflexive,
then there is a
$w^{\ast}$-$\ell sc$ convex $f$
on $X^{\ast}$ and $\Lambda \in X^{\ast}$
such that $f$ is $G$-differentiable at $\Lambda$ and
$f^{\prime}(\Lambda) \in
X^{\ast\ast} \backslash X$ (and, a fortiori, $f$ is not Fr\'echet
differentiable at $\Lambda$).

\smallskip
\pf
Let $Y$ be a separable nonreflexive subspace of $X$.
Let $y^* \in Y^{\ast}$ be such that $y^*$
does not attain its norm on $B_{Y}$.  Let $y^{\ast\ast} \in
S_{Y^{\ast\ast}}$ be
such that $\langle y^{\ast\ast} , y^* \rangle = 1$.
Note that $y^{\ast\ast} \in S_{Y^{\ast\ast}}
\backslash Y$ because $y^*$ does not attain its norm on $B_{Y}$.  By
the
Odell-Rosenthal theorem (see \cite{D}, p.236),
choose $\left\{y_n\right\} \subset
B_Y$ such that $y_n \stackrel{w^{\ast}}{\rightarrow}
y^{\ast\ast}$.  Now $Y^{\ast\ast}
= Y^{\perp\perp} \subset X^{\ast\ast}$
and some careful ``identification checks" show
that $y_n \stackrel{w^{\ast}}{\rightarrow} y^{\ast\ast}$ as elements of
$X^{\ast\ast}$ and $y^{\ast\ast} \in X^{\ast\ast} \backslash X$.
Let $\Lambda$ be a
norm preserving extension of $y^*$, then
$\langle y^{\ast\ast} , \Lambda\rangle = 1$ and we define
$f$ on $X^{\ast}$ by
\[ f(x^*) = \sup \bigl\{ \langle x^*, y_n \bigr\rangle - 1 - a_n : n
\in \Bbb{N}
\bigr\} ~~\mbox{where}~~ a_n \downarrow 0. \]

We now show that $y^{\ast\ast} \in \partial f(\Lambda)$.  Indeed,
\[ \begin{array}{lcl}
f(\Lambda + x^*) - f(\Lambda)= f(\Lambda + x^*)
& = & {\displaystyle \sup^{}_{n}} \bigl\{
\langle \Lambda
+ x^* , y_n \rangle - 1 - a_n \bigr\} \\
& \geq & {\displaystyle \lim^{}_{n \rightarrow \infty}}
\bigl\{ \langle \Lambda , y_n
\rangle - 1 - a_n + \langle x^* , y_n \rangle \bigr\} =
\langle y^{**}, x^* \rangle.
\end{array} \]
To see that $f$ is $G$-differentiable, fix
$x^* \in X^{\ast}$ and let $\epsilon >
0$.  Choose $n_0$ so that $|\langle y^{\ast\ast} - y_n , x^* \rangle |
\leq \epsilon \| x^* \|$ for
$n \geq n_0$.  Now if $2\|t\,x^*\| < \min \left\{ a_1 ,
\ldots , a_{n_{0}} \right\}$, we
have
\[ \begin{array}{lcl}
0 \leq f (\Lambda + t x^* ) - f(\Lambda) -
\langle y^{\ast\ast} , t x^* \rangle & = &
{\displaystyle \sup^{}_{n}} \bigl\{
\langle \Lambda + t x^* ,y_n \rangle - 1 - a_n
\bigr\} - \langle y^{\ast\ast} , t x^* \rangle \\
& = & {\displaystyle \sup_n} \bigl\{ \langle \Lambda, y_n \rangle - 1 +
\langle y_n - y^{**}, tx^* \rangle - a_n \bigr\} \\
& \le & \max\Bigl\{0,{\displaystyle \sup^{}_{n \geq n_0}}
\bigl\{ \langle y_n - y^{**}, t x^* \rangle
- a_n \bigr\}\Bigr\} \\
& \leq & {\displaystyle \sup^{}_{n \geq n_0}}
\bigl\{ | \langle y^{\ast\ast} - y_n ,
t x^* \rangle | \bigr\} \leq \epsilon \|t x^*\|.
\end{array} \]
Thus $f$ is $G$-differentiable at $\Lambda$
with $G$-derivative $y^{\ast\ast} \in
Y^{\ast\ast} \backslash Y$.
\finpf

Using the results of \cite{BF} and \u{S}mulyan's
test type arguments in a fashion similar to Theorem
\ref{thm3.7}, one can also obtain the following result.  We will not
provide the
details.

\begin{thm} \label{thm3.8}
For a Banach space $X$, the following are equivalent.
\begin{itemize}
\item[(a)] $X$ has the Dunford-Pettis property.
\item[(b)] $G$-differentiability and $WH$-differentiability coincide
for
$w^{\ast}$-$\ell sc$, continuous convex functions on $X^{\ast}$.
\item[(c)] $G$-differentiability and $WH$-differentiability coincide
for dual norms
on $X^{\ast}$.
\end{itemize}
\end{thm}

We next consider what happens for ${\cal F}$ a family of norms alone.

\medskip\noindent
{\bf Remark.} \begin{itemize} \item[(a)] $G_{\nm}=F_{\nm}$
on $X$ implies $G_{\dn} =
F_{\dn}$ on $X^{\ast}$, but not conversely.
\item[(b)] $G_{\nm} = WH_{\nm}$ on $X$ implies $G_{\dn} = WH_{\dn}$,
but not conversely.
\item[(c)] $WH_{\nm} = F_{\nm}$ on $X$ does not imply $WH_{\dn} =
F_{\dn}$ on $X^{\ast}$. \end{itemize}

\pf
(a) This is immediate from Theorem \ref{thm3.1}(c) and Theorem 3.7
(since there are Schur spaces that are not finite dimensional).

(b) Since the $DP^{\ast}$ implies the Dunford-Pettis property the first
part
follows from Theorem 3.1(b) and
Theorem \ref{thm3.8}.  However, if $X$ is separable, then by
Corollary \ref{cor3.6}, $H_{\con} = WH_{\con}$
if and only if $X$ is Schur.  Thus the
separable space $C[0,1]$ does not satisfy $G_{\nm} = WH_{\nm}$ yet it
has the Dunford-Pettis property, and thus by Theorem 3.8
satisfies $G_{\dn} = WH_{\dn}$.

(c) On $c_0$ one has $WH_{\nm} = F_{\nm}$ (see Theorem \ref{thm3.1}).
But $\ell_1 = c^{\ast}_{0}$ is a separable dual space and so it admits
a dual
$G$-norm (see \cite {DGZ}, Theorem II.6.7(ii) and
Corollary II.6.9(ii)).  This norm cannot
be everywhere $F$-differentiable since
$\ell_1$ is not reflexive (see \cite{DGZ}, Proposition II.3.4).
However, this dual norm is everywhere
$WH$-differentiable since $\ell_1$ is Schur.
Thus we do not have $WH_{\dn} =
F_{\dn}$ on $\ell_1$.
\finpf

In fact we can be more precise than we were in (c).  Using Theorem 3.7,
Theorem \ref{thm3.8} and Theorem 3.1
along with results from \cite{BF} one can obtain the
following chain of implications.

\centerline{$X$ fails the Schur property but has the
Dunford-Pettis property $\Longrightarrow$}

\centerline{$WH_{\dn} \neq F_{\dn}$ on $X^{\ast}$ $\Longrightarrow$}

\centerline{$X$ fails the Schur
property and $X^{\ast} \supset \ell_1$.}


\end{document}